
\magnification=1200
\def\nline{\par\noindent}

\def\proof{\medskip\noindent{\bf Proof. }}

\def\definition{\noindent{\bf Definition. }}

\def\O{{\cal O}}

\def\V{{\cal V}}

\def\H{{\cal H}}

\def\cN{{\cal N}}
\def\cM{{\cal M}}

\def\ce{{\cal C}}
\def\pe{{\cal P}}

\def\P{{\bf P}}

\def\C{{\bf C}}
\def\M{{\bf M}}
\def\N{{\bf N}}
\def\R{{\bf R}}
\def\T{{\bf T}}

\def\S{\Sigma}

\def\har#1{\smash{\mathop{\hbox to 1 cm{\rightarrowfill}}
\limits^{\scriptstyle#1}_{}}}

\def\dplus#1{\smash{\mathop{\hbox to 0.9 cm{$\quad \oplus$}}
\limits_{\scriptstyle#1}^{}}}

\def\ra{\longrightarrow}
\def\f{\varphi}
\def\iso{\simeq}
\input diagrams


\centerline{{\bf On Euler---Jaczewski sequence}} 
\centerline{{\bf and Remmert---Van de Ven problem for toric varieties}}
\medskip\centerline{{\sl Dedicated to the memory of Krzysztof Jaczewski}}
\medskip\centerline{Gianluca Occhetta and Jaros{\l}aw A. Wi\'sniewski}\par\bigskip\noindent
{\bf Introduction.}\par\medskip\noindent
In [RV] Remmert and Van de Ven posed a problem concerning holomorphic
maps from a complex projective space onto a smooth variety of the same
dimension: the question was whether the target variety has to be the
projective  space as well. The problem has a positive answer provided
by Lazarsfeld, [La].

\proclaim Theorem. {\rm [Lazarsfeld]} Suppose that $Y$ is a smooth
projective variety of positive dimension. If $\f:\P^n\ra Y$ is a
surjective morphism then $Y\iso\P^n$.

Lazarsfeld's proof depends on a (somewhat technical) characterization of
the projective space, by Mori, which was obtained as a by-product of his
proof of Frankel-Hartshorne conjecture, [Mo]. We will explain the result
in the section on rational curves.

Following Lazarsfeld, questions were raised concerning possible extensions
of Remmert -- Van de Ven problem for a broader class of varieties. That
is, given an $X$ from a class of varieties and a morphism $\f: X\ra Y$
onto a smooth projective variety (with possibly additional assumptions on
$\f$ and $Y$, like $\rho(Y)=1$, where $\rho$ denotes the Picard number, or
the rank of Neron-Severi group of the variety), one would like to deduce
the structure of $Y$, preferably to claim that $Y\iso\P^n$, unless $\f$ is
an isomorphism. In particular, the following cases have been considered:  
$X$ is a smooth quadric, see [PS] and [CS], $X$ is an irreducible
symmetric Hermitian space, [Ts], $X$ is rational homogeneous with
$\rho(X)=1$, [HM] and $X$ is a Fano 3-fold, [Am] and [Sc].

%

The main result of the present note is the following. 

\proclaim Theorem 1. Suppose that $X$ is a complete toric variety and $Y$
a smooth projective variety with $\rho(Y)=1$. If $\f: X\ra Y$ is a
surjective morphism then $Y\iso\P^n$.

The proof of the main result follows the line of Lazarsfeld's argument,
in particular we apply Mori's ideas of considering families of rational
curves. A new ingredient of the proof is Euler-Jaczewski sequence for
toric varieties which we explain in the subsequent section. As an
application we obtain a result on varieties admitting two projective
bundle structures (Theorem 2 in the last section of the paper). 

In the present paper all varieties are defined over an algebraically
closed field of characteristic zero.  

\medskip\noindent{\bf Acknowledgements.} The first named author would like
to acknowledge support from Progetto Giovani Ricercatori, Universit\`a di
Milano, which made this collaboration possible. The second named author
would like to acknowledge support from Polish KBN
(2P03A02216).

\par\bigskip\vfill\eject\noindent {\bf Toric varieties.}

\par\medskip\noindent For generalities on toric varities we refer the
reader to [Fu] or [Od]. Let $X=X_\S$ be a toric variety defined
by a fan $\S$ in the space $\N_\R$, with $\N$ denoting the lattice of
1-parameter subgroups of the big torus $\T\iso(\C^*)^n$ and $\M$ denoting
its dual. By $\S^1=\{\rho\}$ let us denote the set of rays (1-dimensional
cones) in the fan $\S$. The equivariant divisor in $X$ (the closure of a
codimension 1 orbit of the action of $\T$) associated to a ray
$\rho\in\S^1$ we shall denote by $D_\rho$. The variety $X$ is decomposed
into the union of the open orbit of $\T$ and the divisor
$\bigcup_{\rho\in\Sigma^1} D_\rho$.

Let us recall the following general fact due to Blanchard, [Bl].

\proclaim Theorem. {\rm [Blanchard]} Let $X$ be a normal complete
variety and $G$ be a connected algebraic group acting on $X$ such that
the induced action on Pic$(X)$ is trivial (assume, for example, that
$H^1(X,\O_X)=0$). Suppose that $\f: X \to Y$ is a morphism to a
projective normal variety with connected fibers, so that
$\f_*\O_X=\O_Y$. Then there exists an action of $G$ on $Y$ such that
$\f$ is equivariant.

\proof \quad Let $L=\f^*L_Y$, where $L_Y$ is an ample line bundle over
$Y$.\nline We have $X=Proj_X(\bigoplus_m L^{\otimes m})$ and
$Y=Proj(\bigoplus_m H^0(X,L^{\otimes m}))$ and $\f: X \to Y$ is induced by
the evaluation morphism $H^0(X,L^{\otimes m})\ra L^{\otimes m}$. The
natural action of $G$ on the graded ring of sections $\bigoplus_m
H^0(X,L^{\otimes m})$ is clearly compatible with the evaluation.

\proclaim Corollary 1. Let $X$ be a toric variety, $\f:X \to Y$ a
morphism to a projective variety $Y$ and  let $X \har{f_0} X' \har{f_1} Y$
be the Stein factorization of $\f$. Then $X'$ is a toric variety.

\proof By Blanchard's theorem the big torus of $X$ acts on $Y$ with an
open orbit. Thus the quotient of the big torus of $X$ by the isotropy of a
general point of $Y$ is a torus and it acts on $Y$ with an open orbit,
hence $Y$ is toric by [Od], Theorem 1.5.\par

\bigskip\noindent{\bf Euler-Jaczewski sequence.}  Let $X$ be a
complete algebraic variety, let $H = H^1(X, \Omega_X)$ and $\H$ be the
sheaf $H \otimes \O_X$.

\par \medskip \definition The short exact sequence $$0 \ra \Omega_X \ra
\pe^{\vee}_X \ra \H \ra 0$$ corresponding to $Id_H \in {\rm Hom}(H,H) =
{\rm Ext}^1(\H,\Omega_X)$ will be called the Euler-Jaczewski sequence of
the variety $X$ whereas the sheaf $\pe_X$ will be called the potential
sheaf of $X$.

\par \medskip\noindent We have the following characterization of toric
varieties in terms of the Euler-Jaczewski sequence, [Ja, 3.1]:

\proclaim Theorem. {\rm [Jaczewski]} A smooth complete and connected
variety $X$ is a toric variety if and only if there exists an effective
divisor $D=\bigcup{i\in I} D_i$ with normal crossing such that $$\pe_X=
\bigoplus_{i \in I}\O_X(D_i)$$ where $D_i$ are the irreducible components
of the divisor $D$. The divisors $D_i$ are then the closures of the
codimension one orbits of the torus.

\noindent The existence of a generalized Euler sequence on a smooth toric variety
was discovered by Batyrev and Melnikov in [BM]. The characterization of
toric varieties in terms of this sequence was proved later by Jaczewski,
[Ja], who apparently was not aware of the Batyrev and Melnikov's work.\par 
\medskip
\noindent On a smooth toric variety we have a short exact sequence: $$0
\ra \M \ra {\rm Div}^TX \ra {\rm Pic}(X) \ra 0$$ where ${\rm Div}^TX$
denotes the group of $\T$ equivariant divisors $\sum_\rho a_\rho
D_\rho$. The first map associates to a character its divisor of poles and
zeroes while the second map associates to a $\T$-equivariant divisor its
class in Pic$(X)$.\nline Dualizing the above sequence and tensoring it
with $\C$ we obtain $$0 \ra N_1(X)_{\C} \ra \bigoplus_{\rho \in \Sigma^1}
\C[\rho] \ra \N_{\C} \ra 0$$ where $N_1(X)_{\C}$ is the space of 1-cycles
on $X$.\nline The maps in the sequence are defined as follows:
$N_1(X)_{\C}\ni Z\ra \sum_{\rho\in\Sigma^1} (Z\cdot D_\rho)\cdot[\rho]$
and $\sum_\rho a_\rho\cdot[\rho]\ra\sum_\rho a_\rho\cdot e_\rho$, where
$e_\rho$ is the generator of the semigroup $\N\cap\rho$.

\proclaim Lemma 1. Let $X=X_\Sigma$ be a toric variety as above.
Consider $\f: X \to Y$ a surjective and generically finite morphism. By
$J(\f)$ let us denote the subset of $\Sigma^1$ corresponding to divisors
$D_\rho$ which are mapped to divisors in $Y$, that is $J(\f) =\{\rho \in
\Sigma^1 |\ \f_*D_\rho \not = 0\}$. Then the restriction of the above map
$\dplus{\rho \in J} \C[\rho] \ra \N_{\C}$ is surjective.

\proof Consider the Stein factorization of $\f:X \to Y$. By corollary 1
the connected part of this factorization is a birational toric morphism
$\f_0:X_{\S} \to X_{\S'}$, where the fan $\S$ is a subdivision of a fan
$\S'$ and the rays of $\S'$ are exactly the rays corresponding to the
divisors $D_\rho$, with $\rho \in J(\f)$. Since the fan $\S'$ is complete,
its rays span $\N_\C$ and we are done.\par \medskip

\noindent Note that, for a toric variety, we have $H^1(X, \Omega_X) =
Pic(X)_{\C}$. In particular $N_1(X)_{\C} \otimes \O_X = \H^{\vee}$;
denoting by $\O_X[\rho]$ the sheaf $\C[\rho]\otimes \O_X$ and by $\cN$ the
sheaf $\N_{\C} \otimes \O_X$ we have a commutative diagram of sheaves
over $X$ with exact rows and columns, which contains both the
Euler-Jaczewski sequence and the
sequence we have just discussed (see [Ja], diagram (7) on page 238):
\newarrow{Isom}=====
\diagram[height=2em,width=2em] &&&&0&&0&&0&&&&\\ &&&&\dTo&&\dTo&&\dTo&&&&\\&& 0
&\rTo&\H^{\vee}&\rTo & \bigoplus \O_X[\rho] & \rTo &\cN & \rTo & 0&&\\&&
&&\dIsom&&\dTo>s&&\dTo>{ev}&&\\&& 0 &\rTo&\H^{\vee}&\rTo &
\bigoplus \O_X(D_\rho) & \rTo & TX & \rTo & 0&&(*)\\&& &&\dTo&&\dTo&&\dTo&&&&\\&&
&& 0 & \rTo & \bigoplus \O_{D_\rho}(D_\rho) & \rIsom & \bigoplus
\O_{D_\rho}(D_\rho) & \rTo & 0&&\\&& &&&&\dTo&&\dTo&&&&\\ &&&&&&0 &&0 &&&&\\ 
\enddiagram \par
\nline Here $s=(s_\rho)$, where $s_\rho$ is the section of
$\O_X(D_\rho)$ vanishing along $D_\rho$, and $ev$ associates to a
1-parameter group $\gamma\in{\bf N}$ its tangent field.

\bigskip\noindent
{\bf Families of curves.}

\medskip\noindent Let $Y$ be a smooth projective variety, and $y \in Y$ a
point. We consider schemes ${\rm Hom}(\P^1,Y)$, parametrizing morphisms
from $\P^1$ to $Y$, and ${\rm Hom}(\P^1,Y; 0 \to y)$, parametrizing
morphisms sending $0 \in \P^1$ to $y \in Y$. Let $V \subset {\rm
Hom}(\P^1,Y)$ be a closed irreducible subvariety; we will call $V$, by
abuse, a family of rational curves on $Y$ and we will denote by $V_y$ the
variety $V \cap {\rm Hom}(\P^1,Y; 0 \to y)$. If the evaluation $F: \P^1
\times {\rm Hom}(\P^1,Y) \to Y$ is a dominating morphism then we will call
$V$ a dominating family of rational curves. Suppose that $\rho(Y)=1$. Then
among all dominating families of rational curves on $Y$ (if there exists
any) we can choose a family $\V$ parametrizing curves of minimal degree
with respect to a chosen ample divisor on $Y$; we will call such a family
a minimal dominating family of rational curves on $Y$.

Let us note that in case of Theorem 1, that is when $Y$ is dominated by
a toric variety and $\rho(Y)=1$, there exists a minimal dominating
family $\V$ of rational curves for $Y$. For this family we will use the 
following version of Mori's result, [Mo]. 

\proclaim Theorem. {\rm [Mori]} Assume that $Y$ is a smooth projective
variety such that $\rho(Y)=1$. Let $\V\subset Hom(\P^1,Y)$ be a minimal
dominating family of rational curves on $Y$. If for a general point $y\in
Y$ and for any $f\in\V_y$ the pull-back $f^*TY$ is ample then $Y\iso
\P^n$.

For the proof of Theorem 1 we will need to understand properties of the
minimal dominating family of rational curves on $Y$ from the point of
view of properties of the dominating variety $X\ra Y$. The key is the
following technical observation which we prove in a more general setup:

\proclaim Lemma 2. Let $\varphi:X \to Y$ be a surjective morphism of
smooth irreducible projective varieties. Assume that $\cM$ is a dominating
family of curves of $Y$, that is, there exists a variety ${\cal C}_Y$ with
morphisms $p: {\cal C}_Y \ra \cM$ and $q: {\cal C}_Y \ra Y$, the latter
morphism is dominating, such that all fibers of $p$ are 1-dimensional and
mapped via $q$ to curves in $Y$. Let ${\cal C}_X$ be an irreducible
component of the fiber product $X \times_Y{\cal C}_Y$ which dominates
${\cal C}_Y$; by $\bar q$ and $\bar\varphi$ let us denote morphisms of
${\cal C}_X$ to $X$ and ${\cal C}_Y$, respectively. Suppose that $D$ is an
irreducible effective Cartier divisor on $X$ such that $\varphi_*D$ is an
ample Cartier divisor on $Y$.  Then, for every $m \in {\cal M}$, the
pull-back $\bar p^*D$ is non zero on every connected component of
$(p\circ\bar\varphi)^{-1}(m)$.

\proof Let  ${\cal C}_X \har{p_0}\cM_0
\har{p_1} {\cal M} $ be the Stein factorization of the map 
$p\circ\bar\f :{\cal C}_X\ra{\cal M}$, so that we the following
commutative diagram
\diagram[height=2em,width=2.5em] &&\cM_0&&\\
&\ruTo^{p_0}&&\rdTo^{p_1}&\\ {\cal C}_X & \rTo_{\bar\varphi} &{\cal
C}_Y & \rTo_p &{\cal M}\\ \dTo<{\bar q}& & \dTo<q &&\\ X& \rTo_\varphi
&Y & & \\ \enddiagram
By our assumptions $p_1 \circ p_0:\bar q^*D \to
\cM$ dominates ${\cal M}$, so its  image in $\cM_0$ is an irreducible
subvariety of maximal dimension, hence, since $\cM_0$ is irreducible, it
coincides with $\cM_0$.

\par\bigskip\noindent{\bf Proof of Theorem 1.}

\par\medskip\noindent First of all, by toric Chow's Lemma we can assume
that $X$ is projective. Then, in view of Corollary 1, considering the
Stein factorization of $\f$, we can assume that $X$ and $Y$ are of the
same dimension. Finally, by taking a desingularization of $X$, we can
assume that $X$ is smooth.

Let $\V$ be a minimal dominating family of rational curves on $Y$.
Let $y \in Y$ be a general point which is not contained in the branch
locus of $\f$ nor it is contained in the image of $\f(\bigcup D_\rho)$.
Let $f:\P^1 \to Y$ be a curve in ${\cal V}$ passing through $y$; the
pullback of the tangent bundle splits into a sum of line bundles:  $f^*TY
\simeq \bigoplus \O_{\P^1}(a_l)$. By Mori's theorem explained above we
shall be done if all numbers $a_l$ are positive.

Suppose, by contradiction that $a_{\bar l} \le 0$ for some $\bar l$; in
this case we have a surjection $f^*TY \ra \O_{\P^1}(a_{\bar l})$ and hence
a surjection $f^*TY \ra \O_{\P^1}$. By Lemma 1 and the commutativity of
diagram $(*)$ the map $\bigoplus_{\rho \in J}\O_X(D_\rho) \ra TX$ is
generically surjective. Since we assume that $\f$ is generically finite,
we get a generically surjective map $\bigoplus_{\rho \in J} \O(D_\rho) \to
\f^*TY$ which, by the choice of $y$, is surjective over $\f^{-1}(y)$.

Let us take ${\cal C}_Y=\V\times\P^1$ with $p:{\cal C}_Y\ra\V$
the projection and $q: {\cal C}_Y\ra Y$ the evaluation.
Now consider the situation discussed in Lemma 2
\diagram[height=2.5em,width=2.5em]
\ce_X &\rTo_{\bar\f} &\ce_Y &\rTo_{p} &\V\\
\dTo<{\bar q}& & \dTo<q&& \\ X & \rTo_\f & Y & &\\ \enddiagram

\nline Let $X_f$ be a connected component of $(p\circ\bar\f)^{-1}(f)$.
We obtain a generically surjective map $\bar q^*(\bigoplus_{\rho \in J}
\O(D_\rho))_{|X_f} \ra ((\bar q^*\circ\f^*)TY)_{|X_f}
=(\bar\f^*(f^*TY))_{|X_f}$ and therefore we have a generically surjective
map $\bar q^*(\bigoplus_{\rho \in J} \O(D_\rho))_{|X_f} \ra \O_{X_f}$.
\nline Thus there exists a non-zero section in 
$$H^0(X_f, \bar q^*(\bigoplus_{\rho \in J}\O(-D_\rho)))=  \bigoplus_{\rho \in J}
H^0(X_f, \bar q^*\O(-D_\rho)),$$
but this is impossible, since $X_f$ is connected and for any
$\rho\in J$ the divisor $\bar q^*D_\rho$ is effective non-zero on $X_f$
by Lemma 2.

\par\bigskip\noindent
{\bf An application.}

\par\medskip\noindent Let $X$ be a smooth variety endowed with two
different $\P$-bundle structures $\f:X \to Y$ and $\psi:X \to Z$. Since
fibers of different extremal ray contractions can meet only in points we
have $\dim X \ge \dim Y + \dim Z$; an easy corollary of Lazarsfeld's
theorem is that we have equality if and only if $X = \P^r \times
\P^s$.\nline Using Theorem 1, we are able to describe the next case; we
have the following

\proclaim Theorem 2. Let $X$ be a smooth projective variety of dimension
$n$, endowed with two different $\P$-bundle structures $\f:X \to Y$ and
$\psi:X \to Z$ such that $\dim Y + \dim Z= n +1$. Then either $n=2m-1$,
$Y=Z=\P^m$ and $X = \P(T \P^m)$ or $Y$ and $Z$ have a $\P$-bundle
structure over a smooth curve $C$ and $X= Y \times_{C} Z$.

\proof Let $[l_{\f}]$ and $[l_{\psi}]$ be the numerical equivalence
classes of lines in the fibers of $\f$ and $\psi$. Let $x$ be a point of
$X$ and let $B_x$ be the set of points of $X$ that can be joined to $x$ by
a chain of rational curves whose numerical class is either $[l_{\f}]$ or
$[l_{\psi}]$.\par \smallskip If there exists a point $x$ such that $B_x=X$
then, by [Ko, IV.3.13.3] the Picard number of $X$ is two, hence $Y$ and
$Z$ are Fano varieties of Picard number one. We also note that every pair
of points of $Z$ is connected by a chain of rational curves whose members
are images of rational curves in $X$ whose numerical class is $[l_{\f}]$.

\nline Let $F_t$ and $F_w$, with $t,\ w\in Z$, be two fibers of $\psi$
such that $\f(F_t) \not = \f(F_w)$. Since points $t$ and $w$ are joined by
a chain of rational curves as above thus there exists a rational curve
$\Gamma$ on $Z$ such that $\f:\psi^{-1}(\Gamma) \ra Y$ is dominant.\nline
Let $\nu: \P^1 \to Z$ be the normalization of $\Gamma$ and let
$X_{\Gamma}\ra\P^1$ be the pull-back of the projective bundle, that is
$X_{\Gamma}= \P^1 \times_{\Gamma} X$. The composition of the induced map
$\bar \nu:X_{\Gamma} \ra X$ with $\f$ is a surjective morphism $\f \circ
\bar \nu: X_{\Gamma} \ra Y$. The variety $X_{\Gamma}$ is a $\P$-bundle on
$\P^1$, hence a toric variety, and $Y$ is a smooth variety with
$\rho(Y)=1$, hence Theorem 1 applies to give $Y \iso \P^{\dim Y}$; in the
same way we also get $Z \iso \P^{\dim Z}$.\nline We conclude this case by
the main
theorem of [Sa], that is we get $\dim Y= \dim Z= m$ and $X \iso
\P(T\P^m)$.

\par \smallskip If $B_x$ is a divisor for every $x \in X$ then, for
general $(x_1, x_2)$ in $X \times X$, the divisors $B_{x_1}$ and $B_{x_2}$
are disjoint and because they are numerically equivalent we have
$B_x^2\equiv 0$. We claim that there exists a regular morphism $p: X\ra C$
with connected fibers, onto a smooth curve $C$, which contracts all
divisors $B_x$. If $H^1(X,\O_X)\ne 0$ then $p$ is obtained from the
Albanese map of $X$ (note that $B_x$ are rationally connected hence
contracted by Albanese). If $H^1(X,\O_X)=0$ then divisors $B_x$ are
linearly equivalent, the linear system $|B_x|$ is base point free and
defines $p$. By construction $p$ factors through $\f$ and $\psi$ to
produce $p_Y: Y\ra C$ and $p_Z: Z\ra C$. Moreover, by construction
$\rho(X/C)=2$ hence $\rho(Y/C)=1$ and $\rho(Z/C)=1$.

\nline
A general $B_x$ is smooth with two projective bundle structures hence, by
what we have said at the beginning of this section, it is a product of
projective spaces. Thus a general fiber of $p_X$ as well as $p_Z$ is a
projective space and over an open Zariski $U$ subset of $C$ both morphisms
are projective bundles. By taking the closure in $Y$ and $Z$ of a
hyperplane section of $p_Y$ and $p_Z$, respectively, defined over the open
set $U$ we get a global relative hyperplane section divisor (we use
$\rho(Z/C)=\rho(Y/C)=1$) hence $p_Y$ and $p_Z$ are projective bundles
globally. The conclusion $X=Y\times_C Z$ is immediate.

\bigskip\noindent
{\bf References.}\par\medskip\noindent

\item{[Am]} Amerik, E., Maps onto certain {F}ano threefolds, {\it Doc.
Math.} {\bf 2} (1997), 195-211.

\item{[Bl]} Blanchard, A. Sur les vari\'et\'es analytiques complexes, {\it
Ann. Sci. \'Ecole Norm. Sup.} {\bf 73}, (1956), 157-202.

\item{[BM]} Batyrev, V. V.; Mel'nikov, D. A.
 A theorem on nonextendability of toric varieties. (Russian) Vestnik Moskov. Univ. Ser. I Mat. Mekh. {\bf 118}
(1986), no. 3, 20--24. English translation: Moscow Univ. Math. Bull. {\bf 41} (1986), no. 3, 23--27.
  
\item{[CS]} Cho, K. and Sato, E., Smooth projective varieties dominated by
smooth quadric hypersurfaces in any characteristic, {\it Math. Z.} {\bf
217}, (1994), 553-565.

\item{[Fu]} Fulton, W., Introduction to toric varieties. Princeton
University Press. Princeton NJ 1993.

\item{[HM]} Hwang, J.-M. and Mok, N., Holomorphic maps from rational
homogeneous spaces of Picard number 1 onto projective manifolds, {\it
Invent. Math.} {\bf 136}, (1999), 208-236.

\item{[Ja]} Jaczewski, K., Generalized Euler sequence and toric
varieties. {\it Classification of algebraic varieties (L'Aquila, 1992)},
227--247, {\it Contemp. Math.} {\bf 162}, AMS 1994

      
\item{[Ko]} Koll\'ar, J., Rational Curves on Algebraic Varieties,
volume~32 of Ergebnisse der Math. Springer Verlag, Berlin, Heidelberg, New
York, Tokio, 1996.

\item{[La]} Lazarsfeld, R., Some applications of the theory of positive
vector bundles. {\it Complete Intersections (Acireale 1983)}, Lecture
Notes in Math., Vol. 1092, Springer, Berlin, (1984), 29-61.

\item{[Mo]} Mori, Sh., Projective manifolds with ample tangent bundles.
{\it Ann. of Math.}  {\bf 110} (1979), 593--606. 

\item{[Od]} Oda, T., Convex bodies and algebraic geometry. An Introduction
to the theory of toric varieties.  volume~15 of Ergebnisse der Math.
Springer Verlag, Berlin, Heidelberg, New York, Tokio, 1988.
    
\item{[PS]} Paranjape, K. H. and Srinivas, V., Self-maps of homogeneous
spaces, {\it Invent. Math.} {\bf 98}, (1989), 425-444.
    
\item{[RV]} Remmert, R. and Van de Ven, A., \"{U}ber holomorphe Abbildung
projektiv-algebraischer Manningfaltigkeiten auf komplexe R\"aume, {\it
Math. Ann.} {\bf 142}, (1961), 453-486.

\item{[Sa]} Sato, E., Varieties which have two projective space bundle
structures, {\it J. Math. Kyoto Univ.} {\bf 25}, (1985), 445-457.

\item{[Sc]} 
Schuhmann, C., Morphisms between Fano threefolds.
{\it J. Algebraic Geom.} {\bf 8}, (1999), 221-244.	
    
\item{[Ts]} Tsai, I H., Rigidity of holomorphic maps between compact
{H}ermitian symmetric spaces, J. Differential Geom. {\bf 33}, (1991),
717-729.

\bigskip
\noindent Gianluca Occhetta \hfill Jaros{\l}aw A.~Wi\'sniewski
\nline
Dipartimento di Matematica ``F. Enriques''\hfill Instytut Matematyki\nline
Universit\`a degli Studi di Milano\hfill Uniwersytet Warszawski \nline
Via Saldini, 50\hfill Banacha 2\nline
I-20133 Milano, Italy\hfill PL-02097 Warszawa, Poland\nline
{\tt occhetta@mat.unimi.it}\hfill{\tt jarekw@mimuw.edu.pl}\nline
\medskip

\end